\documentclass[letterpaper, 10 pt, conference]{IEEEconf}

\usepackage[caption=false]{subfig}
\usepackage[dvips]{graphicx}
\usepackage{amsmath, amssymb}
\usepackage{tikz}
\usetikzlibrary{calc}
\usepackage{psfrag}
\usepackage{pgfplots}

\DeclareMathOperator*{\diag}{diag}
\DeclareMathOperator*{\IQC}{IQC}

\newcounter{thm}
\newcounter{cor}
\newcounter{def}
\newcounter{remcount}
\newtheorem{definition}[def]{Definition}

\newtheorem{theorem}[thm]{Theorem}
\newtheorem{corollary}[cor]{Corollary}
\newtheorem{rem}[remcount]{Remark}

\makeatletter
\renewcommand{\maketag@@@}[1]{\hbox{\m@th\normalsize\normalfont#1}}%
\makeatother

\title{ Distributed Robust Stability Analysis of\\Interconnected Uncertain Systems$^*$}

\author{Martin S. Andersen$^{\dagger}$, Anders Hansson$^{\dagger}$, Sina
  K. Pakazad$^{\dagger}$,  and Anders Rantzer$^{\ddagger}$
\thanks{This work has been supported by the Swedish Department of Education within the ELLIIT project.}
\thanks{Martin
    S. Andersen, Anders Hansson, and Sina K. Pakazad are with the
    Division of Automatic Control, Department of Electrical
    Engineering, Link\"oping University, Sweden. Email:
        {\tt\footnotesize \{martin.andersen, hansson, sina.kh.pa\}@isy.liu.se}. The work was carried out
      while Anders Hansson was a visiting professor at the University of
      California, Los Angeles.}%
\thanks{Anders Rantzer is with the Department of Automatic Control,
    Lund University, Sweden. Email: {\tt\footnotesize anders.rantzer@control.lth.se}}}

\begin{document}
\maketitle
\thispagestyle{empty}
\pagestyle{empty}
\message{The line width is \the\linewidth, and the text width is \the\textwidth}

\begin{abstract}
This paper considers robust stability analysis of a large network of interconnected uncertain systems. To avoid analyzing the entire network as a single large, lumped system, we model the network interconnections with integral quadratic constraints. This approach yields a sparse linear matrix inequality which can be decomposed into a set of smaller, coupled linear matrix inequalities. This allows us to solve the analysis problem efficiently and in a distributed manner. We also show that the decomposed problem is equivalent to the original robustness analysis problem, and hence our method does not introduce additional conservativeness.
\end{abstract}


\section{Introduction}
Robust stability analysis with integral quadratic constraints (IQCs)
provides a unified framework for analysis of uncertain systems with
different kinds of uncertainties \cite{Rantzer}.  The main
computational burden in IQC analysis lies in the solution of a
semi-infinite linear matrix inequality (LMI) which is generally dense.
One method for solving this semi-infinite LMI makes use of the KYP
lemma to reformulate the frequency-dependent semi-infinite LMI as a
single frequency-independent LMI. However, this reformulation yields a
dense LMI, and hence the computational cost is high for large systems
even if the underlying structure is exploited as in \cite{ChapterBook1,Anders,Parrilo,Kao2,and+dah+van10,fuj+kim+koj+oka+yam09,WHJ:09,Wallin,Kao1}. The
semi-infinite LMI can also be solved approximately by discretizing the
frequency variable to obtain a finite number of LMIs. The resulting
LMIs are generally dense, and as a consequence, IQC analysis of
large-scale systems is prohibitively expensive in most cases.

In this paper, we consider IQC-based robustness analysis of a large
network of interconnected uncertain systems. We show that by
representing the network interconnections as quadratic constraints, we
obtain a semi-infinite LMI that is sparse. Reformulating this LMI
using the KYP lemma destroys sparsity, but discretizing the frequency
variable does not, and hence we limit our attention to
frequency-gridding methods. We can exploit sparsity in each LMI in
different ways: (i) we can use an SDP solver that exploits sparsity
\cite{Fukuda_exploitingsparsity,BeY:05,and+dah+van10}, or (ii) we can
use chordal decomposition methods to split a sparse LMI into a set of
smaller, coupled LMIs
\cite{fuj+kim+koj+oka+yam09,kim+koj+mev+yam10}. The first method
solves the robustness analysis problem in a centralized manner, and
hence it requires complete information about all the subsystems in the
network. The decomposition approach, on the other hand, allows us to
solve the robustness analysis problem in a distributed manner where a
small cluster of subsystems only communicates with its neighboring
clusters in the network.

\subsection{Related Work}
Control and robustness analysis of interconnected systems is an active
area of research that has been considered in several papers in the
past few decades; see e.g. 
\cite{P4,Limebeer,P6,Desoer,Langbort}. Different
methods for robustness analysis have been developed, e.g., $\mu$ analysis
and IQC analysis
\cite{RobustAndOptimal,Essentials,Multivariable,Rantzer,UlfIQC}. While
these analysis tools are effective for small and medium-sized
interconnected systems, they fail to produce results for large-scale
interconnected systems because of the high computational cost. To
address this issue, \cite{Kao} and \cite{Ulf} propose an efficient method for
robust stability analysis of interconnected uncertain systems with an
interconnection matrix that is normal, and \cite{Langbort} considers
stability analysis and design methods for networks of certain systems
with uncertain interconnections.  A similar problem is also considered
in \cite{P5}. In \cite{UlfLetter}, the authors consider robust
stability analysis of interconnected uncertain systems using IQC-based
analysis, and they show that when the interconnection matrix is
unitarily diagonalizable, the analysis problem can be decomposed
into smaller problems that are easier to solve. Finally,
\cite{Vinnicombe} shows that by
using Nyquist-like conditions and by considering the dynamics of
individual subsystems and their neighbors, it is possible to relax the
interconnection constraints and arrive at scalable analysis conditions
for interconnected uncertain systems.

\subsection{Outline}
The paper is organized as follows. In
Section~\ref{sec:Interconnection}, we present an integral quadratic
constraint for interconnected uncertain systems, and we show how this
can be used to obtain a sparse analysis problem. In
Section~\ref{sec:Decomposition}, we show how the sparse analysis
problem can be decomposed into smaller, coupled problems, and we
discuss how the analysis problem can be solved distributedly.
Finally, we conclude the paper in Section~\ref{sec:Conclusion}.
\subsection{Notation}
We denote with $\mathbb{R}$ the set of real numbers and with
$\mathbb{R}^{m\times n}$ the set of real $m\times n$ matrices. The
transpose and conjugate transpose of a matrix $G$ are denoted by
$G^{T}$ and $G^{\ast}$, respectively. Given matrices $G^i$, $i =
1,\dots,N$, $\diag(G^1,\dots,G^N)$ denotes a block diagonal matrix
with blocks $G^i$. Similarly, given a set of vectors $v_1,\ldots,v_N$,
the column vector $(v_1,\dots,v_N)$ is obtained by stacking these. The
matrix inequality $G \succ H$ ($G \succeq H$) means that $G-H$ is
positive (semi)definite.  We denote with $\mathcal{L}_2$ the set of
square integrable signals, and $\mathcal{RH}_{\infty}$ represents the
set of real, rational transfer functions with no poles in the closed
right half plane.

\section{Stability Analysis of Interconnected Systems}  \label{sec:Interconnection}
We begin this section with a brief review of some key results from
IQC analysis.

\subsection{IQC Analysis}\label{sec:IQC}
Let $\Delta: \mathbb{R}^d \to \mathbb{R}^d$ be a bounded and causal
operator. The mapping $q = \Delta(p)$ can be characterized using
integral quadratic constraints defined as follow.
\begin{definition}[\cite{Rantzer,UlfIQC}]
Let $\Pi$ be a bounded, self-adjoint operator. The operator
$\Delta$ is said to satisfy the IQC defined by $\Pi$, i.e., $\Delta \in \IQC(\Pi)$, if
\begin{align}
\int_{0}^{\infty} \begin{bmatrix} q \\ \Delta(q) \end{bmatrix}^{T} \Pi \begin{bmatrix} q \\ \Delta(q) \end{bmatrix} \, dt \geq 0, \quad \forall q \in \mathcal{L}_2. 
\end{align} 
This can also be written in frequency domain as 
\begin{align}
\int_{-\infty}^{\infty} \begin{bmatrix} \widehat{q}(j\omega) \\ \widehat{\Delta(q)}(j\omega) \end{bmatrix}^{\ast} \Pi(j\omega) \begin{bmatrix} \widehat{q}(j\omega) \\ \widehat{\Delta(q)}(j\omega) \end{bmatrix} \, d\omega \geq 0, 
\end{align} 
where $\widehat{q}$ and
$\widehat{\Delta(q)}$ are Fourier transforms of the signals $q$ and
$\Delta(q)$, respectively, and $\Pi(j\omega) = \Pi(j\omega)^{\ast}$ is
a transfer function matrix.
\end{definition}

IQCs can be used to describe different classes of operators and hence
different uncertainty sets, e.g., operators with bounded gain, sector
bounded uncertainty, and static nonlinearities \cite{Rantzer}.
We now consider a system described by the equations 
\begin{subequations}\label{eq:UncertainSystem}
\begin{align}\label{eq:sys}
p &= Gq,\\
\label{eq:uncert}
q &= \Delta(p),
\end{align}
\end{subequations}
where $G \in \mathcal{RH}_{\infty}^{m\times m}$ is a transfer
function matrix and $\Delta$ represents the uncertainty in the
system. 
Recall that $\Delta$ is bounded and causal. Robust stability of the system
\eqref{eq:UncertainSystem} can be established using the following
theorem.
\begin{theorem}[\cite{Rantzer,UlfIQC}]\label{thm:thm1}
The uncertain system in \eqref{eq:UncertainSystem} is robustly stable, if 
\begin{enumerate}
\item for all $\tau \in [0,1]$ the interconnection described in \eqref{eq:UncertainSystem}, with $\tau\Delta$, is well-posed.
\item for all $\tau \in [0,1]$, $\tau \Delta \in \IQC(\Pi)$.
\item there exists $\epsilon > 0$ such that
\begin{align}\label{eq:thm1}
\begin{bmatrix} G(j\omega) \\ I \end{bmatrix}^{\ast} \Pi(j\omega) \begin{bmatrix} G(j\omega) \\ I \end{bmatrix} \preceq -\epsilon I, \quad \forall \omega \in \mathbb{R}.
\end{align}
\end{enumerate}
\end{theorem}\vskip1em

It follows from Theorem\ \ref{thm:thm1} that, given a multiplier
$\Pi$, IQC analysis requires the solution of the semi-infinite LMI
\eqref{eq:thm1}. As mentioned in the introduction, we consider an
approximate solution to the IQC analysis problem where the feasibility
of the LMI in \eqref{eq:thm1} is checked only for a finite number of
frequency points. We will see later in this section that this allows us
to reformulate the robust stability problem as a sparse LMI when the
system of interest is a network of uncertain interconnected systems.
\begin{rem}[\cite{UlfIQC,UlfLetter}]\label{rem:DiagonalOperator}
Suppose $\Delta^i \in \IQC(\Pi^i)$ where 
\begin{align}
\Pi^i = \begin{bmatrix} \Pi_{11}^i & \Pi_{12}^i \\ \Pi_{21}^i & \Pi_{22}^i \end{bmatrix}.
\end{align}
The block-diagonal operator $\diag(\Delta^1,\dots,\Delta^N)$ then satisfies the IQC defined by 
\begin{align}\label{eq:IQCDiag}
\tilde{\Pi} = \begin{bmatrix} \bar{\Pi}_{11} & \bar{\Pi}_{12}\\ \bar{\Pi}_{21} & \bar{\Pi}_{22}\end{bmatrix}
\end{align}
where $\bar{\Pi}_{ij} = \diag(\Pi_{ij}^1,\dots,\Pi_{ij}^N)$. 
\end{rem}

\subsection{Network of Uncertain Systems}
\begin{figure}
\begin{center}
\begin{tikzpicture}[scale=0.5, inner sep=0pt]

 \draw (0,0) node [draw, fill=none, minimum size=1.5cm] {};
 \draw (0,3) node (D) [draw, fill=none, minimum size=0.75cm] {$\Delta^i$};

 \draw[->] (-1.5,0.6) node[right=2pt]{$G^i_{pq}$} -- (-2.5,0.6) -- node[left=4pt]{$p^i(t)$} (-2.5,3) -- (D);
 \draw[->] (D) -- (2.5,3) -- node[right=4pt]{$q^i(t)$} (2.5,0.6) -- (1.5,0.6) node[left=2pt]{$G^i_{pw}$};
 \draw[->] (-1.5,-0.6) node[right=2pt]{$G^i_{zq}$} -- (-3,-0.6) node[left=4pt]{$z^i(t)$};
 \draw[->] (3,-0.6) node[right=4pt]{$w^i(t)$} -- (1.5,-0.6) node[left=2pt]{$G^i_{zw}$};

\end{tikzpicture}  
\caption{The $i$th uncertain subsystem with structured uncertainty.}
\label{fig:LFT}
\end{center}
\end{figure}
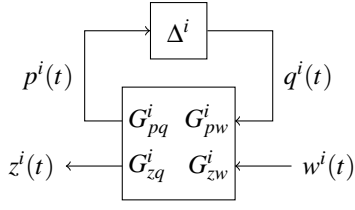

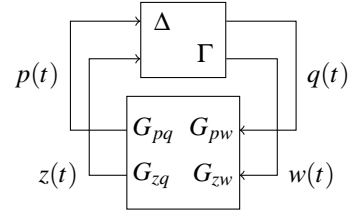
\begin{figure}
\begin{center}
\begin{tikzpicture}[scale=0.5, inner sep=0pt]

 \draw (0,0) node [draw, fill=none, minimum size=1.5cm] {};
 \draw (0,3) node (D) [draw, fill=none, minimum size=1cm] {$\begin{array}{cc} \Delta & \\ & \Gamma \end{array}$};
  
 \draw[->] (-1.5,0.6) node[right=2pt]{$G_{pq}$} -- (-3,0.6) -- node[left=4pt]{$p(t)$} (-3,3.5) -- ($(D.west)+(0,0.5)$);
 \draw[->] ($(D.east)+(0,0.5)$) -- (3,3.5) -- node[right=4pt]{$q(t)$} (3,0.6) -- (1.5,0.6) node[left=2pt]{$G_{pw}$};
 \draw[->] (-1.5,-0.6) node[right=2pt]{$G_{zq}$} -- (-2.5,-0.6) node[left=4pt]{$z(t)$} -- (-2.5,2.5) -- ($(D.west)-(0,0.5)$);
 \draw[->] ($(D.east)-(0,0.5)$) -- (2.5,2.5) -- (2.5,-0.6) node[right=4pt]{$w(t)$} -- (1.5,-0.6) node[left=2pt]{$G_{zw}$};

\end{tikzpicture}
\caption{Compact representation of the network of interconnected systems.}
\label{fig:LFTE}
\end{center}
\end{figure}

Consider a network of $N$ uncertain subsystems of the form
\begin{equation}\label{eq:UncertainIn}
\begin{split}
&p^i = G_{pq}^i q^i + G_{pw}^iw^i \\
&z^i = G_{zq}^i q^i + G_{zw}^iw^i\\
&q^i = \Delta^i(p^i),
\end{split}
\end{equation}
where $G_{pq}^i \in \mathcal{RH}_{\infty}^{d_i \times d_i}$, $G_{pw}^i
\in \mathcal{RH}_{\infty}^{d_i \times m_i}$, $G_{zq}^i \in
\mathcal{RH}_{\infty}^{l_i \times d_i}$, $G_{zw}^i \in
\mathcal{RH}_{\infty}^{l_i \times m_i}$, and $\Delta^i:\mathbb{R}^{d_i} \to \mathbb{R}^{d_i}$. The
$i$th subsystem is shown in Fig.~\ref{fig:LFT}.
The network interconnections are defined by the equation
\begin{align}
  \underbrace{ \begin{bmatrix}
    w^1\\w^2\\ \vdots\\ w^N
  \end{bmatrix}}_{w} = 
\underbrace{  \begin{bmatrix}
    \Gamma_{11} & \Gamma_{12} & \cdots & \Gamma_{1N} \\
    \Gamma_{21} & \Gamma_{22} & \cdots & \Gamma_{2N} \\
    \vdots & \vdots & \ddots & \vdots \\
    \Gamma_{N1} & \Gamma_{N2} & \cdots & \Gamma_{NN} 
  \end{bmatrix}}_{\Gamma}
  \underbrace{ \begin{bmatrix}
    z^1\\z^2\\ \vdots\\ z^N
  \end{bmatrix}}_{z},
\end{align}
where the $ij$th block $\Gamma_{ij}$ is a 0-1 matrix that defines the connections
from system $j$ to system $i$, and $w = (w^1,\dots,w^N)$
and $z = (z^1,\dots,z^N)$ are the stacked inputs and outputs,
respectively. Similarly, we define $q = (q^1,\dots,q^N)$ and $p = (p^1,\dots,p^N)$. Using the interconnection
matrix, the entire system can be expressed as
\begin{subequations}\label{eq:InterSystem}
\begin{align}
p& = G_{pq} q + G_{pw}w \\
z& = G_{zq} q + G_{zw}w\\
q& = \Delta(p)\\
w& = \Gamma z,\label{eq:Interconnection}
\end{align}
\end{subequations}
where the matrices $G_{pq,}$, $G_{pw}$, $G_{zq}$, and $G_{zw}$ are
block-diagonal and of the form $G_{ij} = \diag(G_{ij}^1, \ldots,
G_{ij}^N)$, and $\Delta = \diag(\Delta^1, \dots, \Delta^N)$ is block-diagonal. The full
system is shown in Fig.~\ref{fig:LFTE}.
One approach to robustness analysis of interconnected systems of the
form \eqref{eq:InterSystem} is to eliminate $w$ and $z$ from the system
equations.  This yields the lumped system
\begin{subequations}\label{eq:compact-representation}
\begin{align} 
\label{eq:compact-representation-1}
p& = \bar G q,\\
q& = \Delta(p),
\end{align}
\end{subequations}
where 
\[ \bar G = G_{pq} + G_{pw}(I - \Gamma G_{zw})^{-1}\Gamma G_{zq}\] is
the lumped system matrix. Note that the matrix $I-\Gamma G_{zw}$ must
have a bounded inverse for all frequencies in order for the
interconnection to be well-posed. The lumped system matrix $\bar G$ is
unfortunately dense in general, even if the matrix $I-\Gamma G_{zw}$
is sparse. This follows from the Cayley--Hamilton theorem. As a result,
IQC analysis based on the lumped system description
\eqref{eq:compact-representation} is prohibitively expensive for large
networks of uncertain systems.

\subsection{IQCs for Interconnections}
Robust stability of the interconnected uncertain system \eqref{eq:InterSystem} can be
investigated using the IQC framework. To this end, we now define an
IQC for the interconnection equation $w = \Gamma z$.  The
equation $w = \Gamma z$ can be expressed as the following quadratic
constraint
\begin{align}\label{eq:IQCSimple}
-\|w-\Gamma z \|_X^2 = -\begin{bmatrix} z \\ w \end{bmatrix}^{\ast}
\begin{bmatrix}
  -\Gamma^T \\ I
\end{bmatrix} X
\begin{bmatrix}
  -\Gamma & I
\end{bmatrix}
\begin{bmatrix} z \\ w \end{bmatrix}\geq 0
\end{align}
where $\| \cdot \|_X$ denotes the norm induced by the inner product
$\langle \alpha , X \beta \rangle$ for some $X \succ 0$ of order
$\bar m = \sum_{i=1}^N m_i$.
As a result, we can define an IQC for the interconnections using the following multiplier
\begin{align}
\label{e-iqc-interconnection} 
\hat{\Pi} &=\begin{bmatrix} -\Gamma^T X \Gamma & \Gamma^T X\\ X\Gamma & -X \end{bmatrix}
\end{align}
where $X$ is positive definite.  Note that the framework described in
this paper can be extended to incorporate dynamic and/or uncertain
interconnections by defining suitable IQCs for such interconnections.

We now show that the IQC analysis problem can be expressed as a sparse
LMI by using \eqref{e-iqc-interconnection} to model the
interconnections.  Consider the interconnected uncertain system
defined in~\eqref{eq:InterSystem}, and suppose that
\[ \Delta^i \in \IQC(\Pi^i), \quad i=1,\dots,N\]
such that $\Delta \in \IQC(\bar{\Pi})$ where $\bar{\Pi}$ is
defined in \eqref{eq:IQCDiag}. We will henceforth assume that $\Delta$
satisfies conditions 1 and 2 in Theorem~\ref{thm:thm1}. Then, by
Remark~\ref{rem:DiagonalOperator}, the network of
uncertain systems is robustly stable if there exists $\bar{\Pi}$
and $X \succ 0$ such that
\par{\small\vspace{-1em}
\begin{align} \setlength\arraycolsep{0.15em}
\begin{bmatrix} G_{pq} & G_{pw} \\ G_{zq} & G_{zw} \\ I & 0\\ 0 &
  I \end{bmatrix}^{\ast}
\begin{bmatrix} \bar{\Pi}_{11} &0 & \bar{\Pi}_{12} & 0
  \\ 0 & \hat \Pi_{11} & 0 & \hat{\Pi}_{12}\\ \bar{\Pi}_{21} & 0 & \bar{\Pi}_{22} & 0 \\ 0 & 
\hat{\Pi}_{21}& 0 & \hat{\Pi}_{22}   \end{bmatrix}\begin{bmatrix} G_{pq} & G_{pw} \\ G_{zq} & G_{zw} \\ I & 0\\ 0 & I \end{bmatrix} \preceq -\epsilon I,
\end{align}}%
for $\epsilon > 0$ and for all $\omega \in [0,\infty]$. This can also be written as 
\par{\small\vspace{-1em}
\begin{align}\setlength\arraycolsep{0.15em}
\begin{bmatrix} G_{pq} & G_{pw} \\  I & 0\\G_{zq} & G_{zw} \\ 0 &
  I \end{bmatrix}^{\ast} 
\begin{bmatrix} \bar{\Pi}_{11} &\bar{\Pi}_{12} & 0 & 0 \\ \bar{\Pi}_{21} & \bar{\Pi}_{22} & 0 & 0\\ 0 & 0 & \hat{\Pi}_{11} & \hat{\Pi}_{12} \\ 0 & 
0 & \hat{\Pi}_{21} & \hat{\Pi}_{22}   \end{bmatrix}\begin{bmatrix} G_{pq} & G_{pw} \\  I
& 0\\G_{zq} & G_{zw} \\ 0 & I \end{bmatrix} \preceq -\epsilon I, 
\end{align}}%
or equivalently, as
\begin{multline}\label{eq:IQCInterconnected}
\begin{bmatrix} G_{pq} & G_{pw} \\ I & 0 \end{bmatrix}^{\ast}\begin{bmatrix} \bar{\Pi}_{11} & \bar{\Pi}_{12} \\ \bar{\Pi}_{21} & \bar{\Pi}_{22} \end{bmatrix}\begin{bmatrix} G_{pq} & G_{pw} \\ I & 0 \end{bmatrix} -\\
 \begin{bmatrix} -G_{zq}^{\ast}\Gamma^T\\ I
  -G_{zw}^{\ast}\Gamma^T \end{bmatrix}X\begin{bmatrix} -\Gamma
  G_{zq} & I-\Gamma G_{zw} \end{bmatrix} \preceq -\epsilon I.
\end{multline}
The following theorem establishes the equivalence between robustness
analysis of the lumped system \eqref{eq:compact-representation} via
\eqref{eq:thm1} and robustness analysis of \eqref{eq:InterSystem} via
\eqref{eq:IQCInterconnected}.
\begin{theorem}\label{thm:IQC}
The LMI \eqref{eq:IQCInterconnected} is feasible if and only if 
\begin{align} \label{e-iqc-equiv}
\tilde G_{11} = \begin{bmatrix}
  \bar G\\ I
\end{bmatrix}^\ast \bar \Pi
\begin{bmatrix}
  \bar G \\ I
\end{bmatrix} \preceq -\epsilon I
\end{align}
is feasible.
\end{theorem}
\subsubsection*{proof} We start by applying a congruence transformation $A \mapsto
  T^\ast A T$ to the left-hand side of 
\eqref{eq:IQCInterconnected} where $T$ is nonsingular and defined as
\begin{align*}
T = \begin{bmatrix} I & 0 \\ (I-\Gamma G_{zw})^{-1}\Gamma G_{zq} & I \end{bmatrix}.
\end{align*}
The result is a 2-by-2 block matrix 
\begin{align}\label{eq:IQCTransformed}
\tilde{G} = \begin{bmatrix} \tilde{G}_{11}  &\tilde{G}_{12} \\ \tilde{G}_{21} & \tilde{G}_{22} \end{bmatrix}
\end{align}
where $\tilde{G}_{11}$ is defined as in \eqref{e-iqc-equiv}, and 
\begin{align*}
\Tilde{G}_{12} &= \Tilde{G}_{21}^{\ast} = \begin{bmatrix} \bar{G} \\
  I \end{bmatrix}^{\ast} \bar{\Pi}\begin{bmatrix} G_{pw} \\
  0 \end{bmatrix}, \\
\Tilde{G}_{22} &= \begin{bmatrix}G_{pw} \\ 0\end{bmatrix}^{\ast}
\bar{\Pi}\begin{bmatrix} G_{pw} \\ 0\end{bmatrix} - (I-\Gamma G_{zw})^{\ast}X(I-\Gamma G_{zw}).
\end{align*}
Since $\tilde G$ is homogeneous in $\bar \Pi$ and $X$, it follows that
\eqref{eq:IQCInterconnected} is equivalent to $\tilde G \prec 0$. Now
if \eqref{e-iqc-equiv} is feasible, then \eqref{eq:IQCInterconnected} is feasible if 
\begin{align} \label{e-iqc-schur}
  S = \tilde G_{22} - \tilde G_{21} \tilde G_{11}^{-1} \tilde G_{21}^* \prec 0.
\end{align}
Since $(I-\Gamma G_{zw})$ is nonsingular for any $\bar{\Pi}$ such
that $\tilde{G}_{11} \prec 0$, we can scale $X\succ 0$ so that $S
\prec 0$.
Conversely, if \eqref{eq:IQCInterconnected} is feasible, then $\tilde
G_{11} \prec 0$ and hence we can scale $\bar \Pi$ and $X$ so that
$\tilde G_{11} \preceq -\epsilon I$. $\hspace{5mm} \blacksquare$

Theorem \ref{thm:IQC} implies that our approach does not introduce
conservativeness. However, the LMI \eqref{eq:IQCInterconnected} is
generally dense if the matrix $X$ is dense. We have the following corollary to Theorem
\ref{thm:IQC}.
\begin{corollary} \label{cor:IQC}
  In \eqref{eq:IQCInterconnected}, it is sufficient to consider a 
  diagonal scaling matrix of the form $X= x I$ with $x > 0$.
\end{corollary}
\subsubsection*{proof}
  Suppose $\bar \Pi$ is feasible in \eqref{e-iqc-equiv}. Then we can
  choose $X = x I$ with $x >0$ such that the Schur complement $S$ in
  \eqref{e-iqc-schur} is negative definite.$\hspace{5mm} \blacksquare$

Corollary \ref{cor:IQC} implies that we can choose $X$ to be any
positive diagonal matrix. This mean that the LMI
\eqref{eq:IQCInterconnected} becomes sparse when the interconnection
matrix $\Gamma$ is sufficiently sparse. In the next section, we
discuss how we can solve \eqref{eq:IQCInterconnected} efficiently when
the LMI is sparse.

\section{Decomposing The Analysis Problem}\label{sec:Decomposition}	

In networks where each subsystem is connected only to a small number
of neighboring subsystems, the LMI \eqref{eq:IQCInterconnected} is
generally quite sparse. As mentioned in the introduction, we can solve
the LMI \eqref{eq:IQCInterconnected} in a centralized manner using a
sparse SDP solver such as DSDP, or alternatively, we can make use of
decomposition techniques to facilitate parallel and/or distributed
computation. In this section, we discuss chordal decomposition, and we
show how this decomposition can be used to formulate an efficient distributed
algorithm for robust stability analysis.

\subsection{Chordal Decomposition}
Chordal sparsity plays a fundamental role in many sparse matrix
algorithms \cite{BlP:94}. We say that a matrix is chordal if the
corresponding sparsity graph is chordal, and a graph is called chordal
if all its cycles of length at least four have a chord; see e.g.\
\cite[Ch.~4]{Gol:04}. A clique is a maximal set of vertices that
induce a complete subgraph.  The cliques of the sparsity graph
correspond to the maximal dense principal submatrices in the
matrix. We will use the following result from \cite{AHMR:88} to
decompose the LMI \eqref{eq:IQCInterconnected}.
\begin{theorem}[Agler et al.~\cite{AHMR:88}]\label{e-chordal-decomposition-thm}
  Let $A$ be negative semidefinite and chordal with cliques $J_1,J_2,\ldots,J_L$. Then
  there exists a decomposition
  \begin{align} \label{e-chordal-decomposition}
       A = \sum_{i=1}^L A_i, \quad A_i \preceq 0, \ i =1,\ldots,L
  \end{align}
  such that $[A_i]_{jk} = 0$ for all $(j,k) \notin J_i \times J_i $.
\end{theorem}
\subsubsection*{proof}
  See \cite{AHMR:88} and \cite{Kak:10}.$\hspace{5mm} \blacksquare$

\begin{rem}
  It is also possible to apply Theorem
  \ref{e-chordal-decomposition-thm} to negative semidefinite matrices
  with nonchordal sparsity by using a so-called chordal embedding. A
  chordal embedding can be computed efficiently with symbolic
  factorization techniques; see e.g. \cite{OCF:76}.
\end{rem}

A band matrix is an example of a chordal matrix, i.e., if $A$ has
half-bandwidth $w$, then
\[ A_{jk}  = 0, \quad |j-k| > w, \]
and the cliques of the corresponding sparsity graph are given by the sets
\[ J_i = \{i,i+1,\ldots, i+w\}, \quad i = 1,\ldots,n-w, \] where $n$
is the order of $A$. 

The cliques of a chordal graph can be represented
using a clique tree which is a maximum-weight spanning tree of a
weighted clique intersection graph; refer to \cite{BlP:94} for further
details. The clique tree can be used to parameterize all possible
clique-based splittings of $A$. For example, consider an irreducible
chordal matrix $A\preceq 0$ with two cliques $J_1$ and $J_2$, and let
$A = \tilde{A}_1 + \tilde{A}_2$ be an any clique-based splitting of
$A$ such that $[\tilde{A}_i]_{jk} = 0$ for all $(j,k) \notin J_i
\times J_i $. Then, by Theorem~\ref{e-chordal-decomposition-thm}, there exists a symmetric matrix $Z$ that satisfies
\[ Z_{ij} = 0, \quad (i,j) \notin (J_1 \cap J_2) \times (J_1 \cap J_2)\]
such that
\[ A_1 = \tilde{A_1} + Z \preceq 0, \ A_2 = \tilde A_2 - Z \preceq
0.\] It is easy to verify that the matrices $A_1$ and $A_2$ satisfy
\eqref{e-chordal-decomposition}.
We can use this decomposition technique to parameterize all possible
splittings of the LMI \eqref{eq:IQCInterconnected}. Suppose
\eqref{eq:IQCInterconnected} is irreducible and chordal with cliques $J_1,\ldots, J_L$.
Then if we choose $X = x I$, we can split
\eqref{eq:IQCInterconnected} as 
\begin{align}\label{e-chordal-decomposition-disjoint}
\tilde A_1+\tilde A_2+\cdots+\tilde A_L  
\end{align}
 such that the matrices $\tilde A_1,\ldots,\tilde A_L$ share
only the variable $x$, and $[\tilde A_i]_{jk} = 0$ for all $(j,k) \notin
J_i \times J_i $. This means that \eqref{eq:IQCInterconnected} is
equivalent to $L$ coupled LMIs
\begin{align}\label{e-chordal-decomposition-coupling}
  \tilde{A}_i + Z_{p_ii} - \sum_{j \in \mathrm{ch}(i)} Z_{ij} \preceq
  0, \ i=1,\ldots,L
\end{align}
where $p_i$ is the parent of the $i$th clique in the clique tree and 
$\mathrm{ch}(i)$ is the set of children of clique $i$. The matrix $Z_{ij}$
couples cliques $i$ and $j$, and it satisfies
\[ [Z_{ij}]_{kl} = 0, \quad (k,l) \notin (J_i \cap J_j) \times (J_i
\cap J_j). \] Note that the term $Z_{p_i i }$ in
\eqref{e-chordal-decomposition-coupling} disappears if clique $i$ is
the root of the clique tree, and the summation over the children of
clique $i$ disappears if $i$ is a leaf clique. Furthermore, note that
the sum of the left-hand sides of the LMIs in
\eqref{e-chordal-decomposition-coupling} is equal to
\eqref{e-chordal-decomposition-disjoint}.


\subsection{Distributed Algorithm}
The set of LMIs \eqref{e-chordal-decomposition-coupling} can be 
expressed as
\begin{subequations}\label{e-decomp-consensus}
\begin{gather}\label{e-decomp-consensus-1}
    \tilde{A}_i + \hat Z_{p_ii} - \sum_{j \in \mathrm{ch}(i)} Z_{ij} \preceq
  0, \ i=1,\ldots,L \\
  \label{e-decomp-consensus-2}
  x_k = x_l, \ \hat Z_{kl} = Z_{kl}, \ \forall (k,l) \in \mathcal T
\end{gather}
\end{subequations}
where $\hat Z_{p_i i}$ is the copy of $Z_{p_i i}$ associated with
clique $i$, $x_i$ is the copy of $x$ associated with clique $i$, and $(k,l)
\in \mathcal T$ means that $(k,l)$ is an edge in the clique tree
$\mathcal T$. The coupling between the LMIs
\eqref{e-decomp-consensus-1} is now described as a set of 
equality constraints
\eqref{e-decomp-consensus-2}. Fig.~\ref{fig:subtree} shows a subtree
of the clique tree with the auxiliary variables associated with the subtree. 
\begin{figure}
  \centering
  \begin{tikzpicture}
	\tikzstyle{nodeA} = [draw, circle, fill=none, minimum size=6mm]
	\draw node[style=nodeA] (p) at (0,1.5) {\small $p$};
	\draw node[style=nodeA, thick] (i) at (0,0) {\small $i$};  
	\draw node[style=nodeA] (j1) at (-2,-1.5) {\small $c_1$};  
	\draw node[style=nodeA] (j2) at (-0.6,-1.5) {\small $c_2$};  
	\draw node[style=nodeA] (jk) at (2,-1.5) {\small $c_k$};  
	\draw (p) -- node[right=2pt,near start]{\scriptsize $Z_{pi}$} (i) node[right=2pt,near end]{\scriptsize $\hat Z_{pi}$};
	\draw (i) -- node[above left=-2pt,near start]{\scriptsize $Z_{ic_1}$} (j1) node[above left=-2pt,near end]{\scriptsize $\hat Z_{ic_1}$};
	\draw (i) -- node[right=0pt,near start]{\scriptsize $Z_{ic_2}$} (j2) node[right=0pt,near end]{\scriptsize $\hat Z_{ic_2}$};
	\draw (i) -- node[above right=-2pt,near start]{\scriptsize $Z_{ic_k}$} (jk) node[above right=-2pt,near end]{\scriptsize $\hat Z_{ic_k}$};
	
	\draw node at (0.7,-1.5) {$\cdots$};
\end{tikzpicture}
  \caption{Subtree of clique tree: clique $i$, its parent, and its children. }
  \label{fig:subtree}
\end{figure}
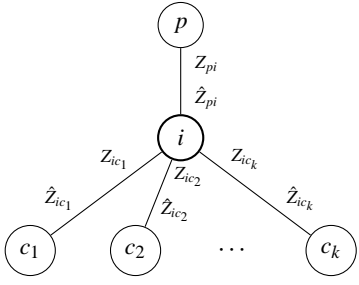
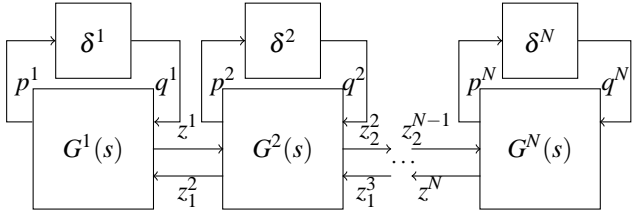
\begin{figure}[t]
\begin{center}
\begin{tikzpicture}[scale=0.36, inner sep=0pt]

  \tikzstyle{nodeA} = [draw, fill=none, minimum size=1.6cm]
  \tikzstyle{nodeB} = [draw, fill=none, minimum size=1.0cm]

  \draw (0,0) node [style=nodeA]
  {$G^1(s)$};
  \draw (7,0) node [style=nodeA]
  {$G^2(s)$};
  \draw (16.5,0) node [style=nodeA]
  {$G^N(s)$};

  \draw (0,4)  node [style=nodeB]
  {$\delta^1$};
  \draw (7,4)  node [style=nodeB]
  {$\delta^2$};
  \draw (16.5,4) node [style=nodeB]
  {$\delta^N$};
 
  \draw[->] (-2.25,1) -- (-3.2,1) -- node[right]{$\phantom{.}p^1$} (-3.2,4) -- (-1.4,4);
  \draw[->] (1.4,4) -- (3.2,4) -- node[left]{$q^1$} (3.2,1) -- (2.25,1);
  \draw[->] (2.2,0) -- node[above]{\raisebox{1.5mm}{$z^1$}} (4.8,0);
  \draw[->] (4.8,-1) -- node[below]{\raisebox{-2.5mm}{$z_1^2$}} (2.2,-1);

  \draw[->] (4.8,1) -- (4,1) -- node[right]{$\phantom{.}p^2$} (4,4) -- (5.6,4);
  \draw[->] (8.4,4) -- (10.1,4) -- node[left]{$\ q^2$} (10.1,1) -- (9.2,1);
  \draw[->] (9.2,0) -- node[above]{\raisebox{1.5mm}{$\ z_2^2$}} (11,0);
  \draw[->] (11,-1) --  node[below]{\raisebox{-2.5mm}{$z_1^3$}} (9.2,-1);

  \draw[->] (14.3,1) -- (13.5,1) -- node[right]{$\phantom{.}p^N$} (13.5,4) -- (15.1,4);
  \draw[->] (17.9,4) -- (19.85,4) -- node[left]{$\ q^N$} (19.85,1) -- (18.7,1);
  \draw[->] (11.75,0) -- node[above]{\raisebox{1.5mm}{$z_2^{N-1}\hspace{5mm}$}} (14.3,0);
  \draw[->] (14.3,-1) -- node[below]{\raisebox{-2.5mm}{$z^N\hspace{4mm}$}} (11.75,-1);

  \draw (11.5,-0.5) node [fill=none] {{\small  $\cdots$}};

\end{tikzpicture} 
\caption{A chain of $N$ uncertain system.}
\label{fig:System}
\end{center}
\end{figure}
The formulation \eqref{e-decomp-consensus} is a convex feasibility
problem of the form
\begin{align} \label{e-conv-feas}
  \begin{array}{ll}
    \mbox{find} & x, s_1,s_2,\ldots,s_L   \\
    \mbox{subject to} & s_i \in \mathcal C_i, \ i=1,\ldots, L \\
                                  & s_i = H_i(z), \ i=1,\ldots, L
  \end{array} 
\end{align}
where $\mathcal C_1,\ldots,\mathcal C_L$ are convex sets, $s_1,\ldots,
s_L$ are local variables, and the constraints $s_i = H_i(z)$ ensure
global consensus. The problem \eqref{e-conv-feas} can be solved
distributedly using
e.g.~the alternating projection (AP) method, Dykstra's AP method, or
the alternating direction method of multipliers \cite{BoydADMM}.

\subsection{Example: Chain of Uncertain Systems}\label{sec:Chain}
We now consider as an example a chain of $N$ uncertain systems where
each system $G^i(s)$ is defined as in \eqref{eq:UncertainIn} with
scalar uncertainties $\delta^1,\ldots,\delta^N$.
Fig.~\ref{fig:System} shows the chain of uncertain systems.

\begin{figure}
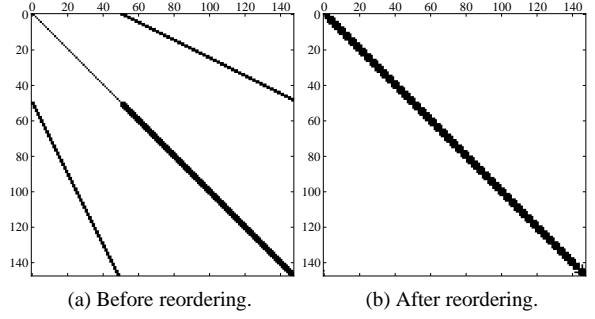

  \centering
  \subfloat[Before reordering.]{\input{band_spy_50.tex}}
  \subfloat[After reordering.]{\input{band_spy_50_perm.tex}}
  \caption{Sparsity pattern associated with \eqref{eq:IQCInterconnected} for a chain of 50 systems.}
  \label{fig:band_sparsity}
\end{figure}

The inputs and outputs are defined as $w_i = ( w_i^1 , w_i^2 )$ and
$z_i = ( z_i^1 , z_i^2 )$ for $i=2,\ldots,N-1$, and moreover, $w_1
= w_1^2$, $w_N = w_{N}^1$, $z_1 = z_1^2$, and $z_N = z_N^1$. The
interconnections are described by the equations $w_i^1 = z_{i-1}^2$
and $w_i^2 = z_{i+1}^1$ for $i=2,\ldots,N-1$, and $w_1^2 = z_2^1$ and
$w_N^1 = z_{N-1}^2$. We assume that the signals $w_i^j$ and $z_i^j$
are scalar-valued. As a result, the interconnection matrix $\Gamma$
that describes this network has nonzero blocks $\Gamma_{i,i-1} =
\Gamma_{i-1,i}^T$ for $i = 2,\ldots,N$, and these blocks are given by
\begin{align}
  \Gamma_{i,i-1} = \Gamma_{i-1,i}^T =
  \begin{bmatrix}
    0 & 1 \\ 0 & 0
  \end{bmatrix}, \ i=3,\ldots,N-1,
\end{align}
and
\begin{align}
  \Gamma_{21} = \Gamma_{12}^T = (1,0), \quad \Gamma_{N-1,N} =
  \Gamma_{N,N-1}^T = (0,1).
\end{align}
In this example, we will assume that the uncertainties
$\delta^1,\ldots,\delta^N$ are scalar, and moreover, $\delta^i \in
\IQC(\Pi_i)$ for all $i =1,\ldots,N$. If we let $X$ be diagonal, the
sparsity pattern associated with the LMI \eqref{eq:IQCInterconnected}
is then chordal with $2N-2$ cliques, and the largest clique is of
order 4. The LMI \eqref{eq:IQCInterconnected} can therefore be
decomposed into $2N-2$ coupled LMIs of order at most 4. The sparsity
pattern associated with \eqref{eq:IQCInterconnected} for a chain of 50
uncertain systems is shown in Fig.~\ref{fig:band_sparsity}. In this
example, it is also possible to combine overlapping cliques such that
we get a total of $N$ cliques of order at most $5$, and in general,
there is a trade-off between the number of cliques and the order of
the cliques.
Note that analyzing the lumped system yields an LMI
\eqref{e-iqc-equiv} of order $N$ whereas the sparse LMI \eqref{eq:IQCInterconnected}
is of order $3N-2$. For large networks, solving the sparse LMI can be much faster, but for small and medium- sized networks, the original dense LMI  \eqref{e-iqc-equiv} may be cheaper to solve.

\section{Conclusions}\label{sec:Conclusion}
IQC-based robustness analysis of a large network of interconnected
systems involves the solution of a large, dense LMI. By expressing the
network interconnections in terms of IQCs, we have shown that it is
possible to express the robustness analysis problem as a sparse LMI
that can be decomposed into a set of smaller but coupled LMIs. The
decomposed problem can be solved using distributed computations, and
we have shown that it is equivalent to the original problem.  These
findings suggest that our method is applicable to robustness analysis
of large networks of interconnected uncertain systems, but further
work needs to be done to establish what types of network structure
yield computationally efficient decompositions.

\bibliographystyle{IEEEtran}  
\bibliography{IEEEabrv,CDCPaper} 

\end{document}